\documentclass[11pt]{article}

\usepackage{titling}
\usepackage{url,color,tikz}
\usepackage{mathrsfs}
\usepackage{pgfplots}
\usepackage{graphicx}
\usepackage{subcaption}
\usepackage{bbm}
\usepgfplotslibrary{groupplots}
\usepackage{tikzscale}

\usepackage[english]{babel}
\usepackage[utf8]{inputenc}

\usepackage{amsmath,amsthm,amssymb}
\usepackage{graphicx}
\usepackage[T1]{fontenc}
\usepackage[a4paper,bindingoffset=0.5cm,left=2cm,right=2cm,top=2cm,bottom=2cm,footskip=.8cm]{geometry}
\usepackage[md]{titlesec}
\usepackage{color}

\usepackage{mathrsfs}  
\usepackage[mathscr]{euscript}

\usepackage{mathrsfs}  

\definecolor{light}{gray}{0.8}
\def\JGB#1{\textcolor{red}{#1}}
\def\ERASE#1{\textcolor{light}{#1}}

\usepackage[english]{babel}
\usepackage{amsmath,amsthm,amssymb}
\usepackage{graphicx}
\usepackage[a4paper,bindingoffset=0.5cm,left=2cm,right=2cm,top=2cm,bottom=2cm,footskip=.8cm]{geometry}
\usepackage{color,tikz}

\usepackage{bm}
\newcommand{\bs}{\pmb}

\newtheorem{thm}{Theorem}

\newtheorem{cor}{Corollary}
\newtheorem{prop}{Proposition}

\theoremstyle{remark}
\newtheorem{rem}{Remark}
\newtheorem{exmp}{Example}

\newcommand{\vt}{\vartheta}

\newcommand{\vb}{\vspace{3.2mm}}
\newcommand{\s}{^\star}

\newcommand{\ami}{{(-,\mbox{\scriptsize {\sc i}})}}
\newcommand{\amii}{{(-,\mbox{\scriptsize{\sc ii}})}}
\newcommand{\api}{{(+,\mbox{\scriptsize{\sc i}})}}
\newcommand{\apii}{{(+,\mbox{\scriptsize{\sc ii}})}}

\title{Refined large deviations asymptotics\\for Markov-modulated infinite-server systems}

\author{Joke\ Blom$\,^\star$, Koen\ De Turck$\,^\dagger$, Michel Mandjes$\,^{\bullet,\star}$}

\date{\today}

\begin{document}
%

\begin{titlepage}
    \maketitle
\begin{itemize}
\item[$^\bullet$] Korteweg-de Vries Institute for Mathematics,
University of Amsterdam, Science Park 904, 1098 XH Amsterdam, the Netherlands.
\item[$^\star$] CWI, P.O. Box 94079, 1090 GB Amsterdam, the Netherlands.
\item[$^{\dagger}$] {Laboratoire Signaux et Syst\`emes (L2S, CNRS UMR8506), \'Ecole CentraleSup\'elec, Universit\'e Paris Saclay, 3 Rue Joliot Curie, Plateau de Moulon, 91190 Gif-sur-Yvette, France.}
\item[email:] Joke.Blom@cwi.nl; koen.deturck@centralesupelec.fr; m.r.h.mandjes@uva.nl
\item[tel:] +33 (0)1 69 85 14 63
\end{itemize}

\end{titlepage}

\maketitle

\begin{abstract}\noindent
Many networking-related settings can be modeled by Markov-modulated infinite-server systems. In such models, the customers' arrival rates and service rates are modulated by a Markovian background process; additionally, there are infinitely many servers (and consequently the resulting model is often used as a proxy for the corresponding many-server model).
The Markov-modulated infinite-server model hardly allows any explicit analysis, apart from results in terms of systems of (ordinary or partial) differential equations for the underlying probability generating functions, and recursions to obtain all moments. As a consequence, recent research efforts have pursued an asymptotic analysis in various limiting regimes, notably the central-limit regime (describing fluctuations around the average behavior) and the large-deviations regime (focusing on rare events). Many of these results use the property that the number of customers in the system obeys a Poisson distribution with a random parameter. 

\noindent The objective of this paper is to develop techniques to accurately approximate tail probabilities in the large-deviations regime. We consider the scaling in which the arrival rates are inflated by a factor $N$, and we are interested in the probability that the number of customers exceeds a given level $Na.$
Where earlier contributions focused on so-called {\it logarithmic asymptotics} of this exceedance probability (which are inherently imprecise), the present paper improves upon those results in that {\it exact asymptotics} are established. These are found in two steps: first the distribution of the random parameter of the Poisson distribution is characterized, and then this knowledge is used to identify the exact asymptotics. The paper is concluded by a set of numerical experiments, in which the accuracy of the asymptotic results is assessed.

\begin{itemize}
\item[$^\bullet$] Korteweg-de Vries Institute for Mathematics,
University of Amsterdam, Science Park 904, 1098 XH Amsterdam, the Netherlands. 
\item[$^\star$] CWI, P.O. Box 94079, 1090 GB Amsterdam, the Netherlands.
\item[$^{\dagger}$] {Laboratoire Signaux et Syst\`emes (L2S, CNRS UMR8506), \'Ecole CentraleSup\'elec, Universit\'e Paris Saclay, 3 Rue Joliot Curie, Plateau de Moulon, 91190 Gif-sur-Yvette, France.}
\end{itemize}

\noindent
M.\ Mandjes is also with  E{\sc urandom}, Eindhoven University of Technology, Eindhoven, the Netherlands, 
and IBIS, Faculty of Economics and Business, University of Amsterdam,
Amsterdam, the Netherlands. M. Mandjes' research is partly funded by the NWO Gravitation project NETWORKS, grant number 024.002.003.

\vb


\end{abstract}


%

\section{Introduction, notation, and preliminaries}
Consider an infinite-server queue modulated by a finite-state irreducible continuous-time Markov chain $J$: when the so-called {\it background process} $J$ is in state $i\in\{1,\ldots,d\}$, jobs arrive according to a Poisson process with rate $\lambda_i$, while the departure rate is $\mu_i$.
These Markov-modulated infinite-server queues have attracted some attention during the past decades; see e.g.\ the early contributions \cite{DAURIA, KEILSONSERVI1993, OCINNEIDEPURDUE} and later \cite{FRALIXADAN}. Importantly, considerably fewer results are available for this model than for the corresponding {\it single}-server queue.
This is primarily due to the fact that, despite the system's simple structure, the Markov-modulated infinite-server queue hardly allows any explicit analysis: whereas the Markov-modulated single-server queue has a matrix-geometric stationary distribution, no such result applies to its infinite-server counterpart. 
The results obtained so far are implicit, in that they are in terms of partial differential equations characterizing the probability generating functions related to the system's transient behavior, and recursions for the corresponding moments (where in each step of the recursion a system of non-homogeneous ordinary differential equations needs to be solved). 

The Markov-modulated infinite-server queue can be applied in various domains, ranging from biology to the performance analysis of particular communication networks. In the present paper the focus lies on the latter application, where the model with an infinite number of servers typically serves as a proxy for its counterpart with a large but finite number of servers. The Markov modulation of the arrival rates and service rates facilitates the modeling of some sort of ?burstiness?; although the concept of Markov modulation has been around for a few decades, it still spurs a considerable amount of research effort \cite{GH, MOR}. For instance, the model can be used to describe the fluctuations in the users' activity level (where each user alternates between transmitting data or being silent). Also, e.g.\ in a wireless setting, the modulation of the service rate can represent  channel conditions that vary over time. In the context of communication networks, a particularly relevant feature concerns {\it rare events}. More specifically, a high activity level corresponds to congestion, and therefore the system should be designed such that such high activity levels occur relatively infrequently.

\vb

Given that, as argued above, explicit analysis is hardly possible, recent research efforts have focused on the exploration of various limiting regimes.
In the first place, significant progress has been made in terms of the derivation of (functional) central limit theorems under specific parameter scalings. When inflating the arrival rates by a factor $N$, and speeding up the background process by a factor $N^\alpha$ (for some $\alpha>0$), in e.g.\  \cite{DAVE,PEIS,MMOR} it has been proven that  the (transient as well as stationary) number of jobs present in the system is, after centering and normalizing, asymptotically Normally distributed. An interesting dichotomy was identified, in that the regimes $\alpha<1$ and $\alpha>1$ lead to qualitatively different asymptotics.

Also the large-deviations regime has been explored, resulting in so-called {\it logarithmic asymptotics} 
{\cite{BKMT,BM,BTM}}. 
In these papers the arrival rates are scaled by a factor $N$ {and} the background process is
{either }left unchanged {or accelerated by a factor $N^{1+\varepsilon},\;\varepsilon>0$}. 
With $M^{(N)}(t)$ the number of jobs present at time $t$ in the resulting system, these papers determine the limit
\begin{equation}
\label{logas}\lim_{N\to\infty}\frac{1}{N}\log p_t^{(N)}(a)=:-I(a),\:\:\mbox{with}\:\:
p_t^{(N)}(a):={\mathbb P}\left(M^{(N)}(t)\ge Na\right),\end{equation}
as well as the corresponding limit for $M^{(N)}(t)$'s steady-state counterpart $M^{(N)}$.
It is observed that these asymptotics are inherently imprecise, as they essentially just entail that
\[p_t^{(N)}(a) = e^{-NI(a)} \Psi(N),\]
for some {\it unknown} subexponential function $\Psi(N)$; we only know that $\Psi(N)$ has the property that, as $N\to\infty$,
\begin{equation}
\frac{1}{N}\log \Psi(N)\to 0.\label{subex}
\end{equation}
Observe that (\ref{subex}) still leaves a substantial amount of freedom: $\Psi(N)$ could be for instance a constant, but also any polynomial function of $N$, or even `big functions' of the type 
$10^6\cdot \exp({N^{0.99}}).$ We conclude that logarithmic asymptotics of the type (\ref{logas}) typically 
provide  valuable insight into the system's rare-event behavior, but that they may be to{o} inaccurate to be used for performance evaluation purposes. This shows that there is a clear need for more precise asymptotic results.

\vb

The main contribution of the present paper is to improve the logarithmic asymptotics (\ref{logas}) to so-called {\it exact} asymptotics: we identify an explicit function $\zeta(\cdot)$ such that, as $N\to\infty$,
\[\frac{p_t^{(N)}(a)}{\zeta(N)}\to 1.\]
As it turns out, this $\zeta(N)$ is the product of the exponential term identified above ($e^{-NI(a)}$), a polynomial term (which is typically of the form $N^{-C}$, for some $C>0$), and a constant. The proof of this property consists of two steps, and relies on the property that $M^{(N)}(t)$ obeys a Poisson distribution with random parameter (as was observed in e.g.\ \cite{BKMT,DAURIA}). 
\begin{itemize}
\item[$\circ$]
In the first step a system of partial differential equations is set up for the distribution of this Poisson parameter. 
\item[$\circ$] In the second step, this is combined with (a uniform version) of the classical result by Bahadur and Rao \cite{BR,HOG} on the exact tail asymptotics of sample means of i.i.d.\ random variables, so as to obtain the exact asymptotics of the tail probability of our interest. 
\end{itemize}

\vb

{\it Model and notation.} As mentioned above, $\lambda_i$ is the (Poissonian) arrival rate when the background process is in state $i$. We let 
\[Q=(q_{ij})_{i,j=1}^d\]
be the ($d {\times} d$) transition rate matrix of the (irreducible) background process $J$, with ${\bs \pi}$ denoting the corresponding invariant probability measure (which is a $d$-dimensional vector ${\bs \pi}$). The entries of $Q$ are non-negative, except for those on the diagonal; the row-sums are assumed to be $0$, where we define $q_i:=-q_{ii}\ge 0$.

Concerning the departure process, two models are considered. In the first, referred to  as Model~{\sc i}, each job present is experiencing a departure rate $\mu_i$ when $J$ is in state $i$; as a consequence, this hazard rate may change during the job's sojourn time (that is, when the background process makes a transition). In the second, Model~{\sc ii}, the crucial difference is that the job's sojourn time is sampled upon arrival: when the background process is then in state $i$, it has an exponential distribution with mean $1/\mu_i$. The evident independence assumptions are imposed.

\vb

{\it Preliminaries.} 
In Model {\sc i}  and {\sc ii}, we have that $M^{(N)}(t)$ has a mixed Poisson distribution, i.e., a Poisson distribution with random parameter \cite{BKMT,DAURIA}. More specifically, with $P(b)$ denoting a Poisson random variable with mean $b>0$, our target probability $p^{(N)}_t(a)$ equals the probability 
${\mathbb P}(P(N\phi_t(J)) \ge Na)$ in Model~{\sc i} and 
${\mathbb P}(P(N\psi_t(J)) \ge Na)$ in Model~{\sc ii},
where the functionals $\phi_t(J)$ and $\psi_t(J)$ of the path
$J\equiv \{J(s): s\in[0,t]\}$
are given by, respectively,
\[\phi_t(J) := \int_0^t \lambda_{J(s)} e^{-\int_s^t \mu_{J( r)} {\rm d}r} {\rm d}s\:\:\:\mbox{
and}
\:\:\:\psi_t(J) := \int_0^t \lambda_{J(s)} e^{-(t-s) \mu_{J(s)}} {\rm d}s.\]
An intuitive explanation for this property is the following. In Model {\sc ii} the probability of a job that has arrived at time $s$ is still present at time $t\in(s,\infty)$ is
\[e^{-(t-s)\mu_{J(s)}},\]
as $\mu_{J(s)}$ is its hazard rate during its entire lifetime. In Model {\sc i} this hazard rate may change over time, in the sense that when the background process is in state $i$ it is $\mu_i$; therefore,
the probability of a job that has arrived at time $s$ is still present at $t$ is
\[e^{-\int_s^t \mu_{J( r)} {\rm d}r} .\]In an earlier paper \cite{BKMT} we have developed a technique to determine for Model {\sc i} numbers $a_t^\ami$ and $a_t^\api$ (such that $0\le a_t^\ami\le a_t^\api$)
being the smallest, resp.\ largest numbers that $\phi_t(J)$ can attain. The analogous result for $\psi_t(J)$ (featuring in Model {\sc ii}) has been presented in  \cite{BM}, resulting in numbers
$a_t^\amii$ and $a_t^\apii$. 

In Model~{\sc ii}, the bounds $a_t^\amii$ and $a_t^\apii$ are explicitly given:
\begin{equation}\label{MAXP} a_t^\amii = \int_0^t \left(\min_{i\in\{1,\ldots,d\}} \lambda_i e^{-(t-s)\mu_i}\right) {\rm d}s,\:\:\:\: a_t^\apii = \int_0^t \left(\max_{i\in\{1,\ldots,d\}} \lambda_i e^{-(t-s)\mu_i} \right){\rm d}s.\end{equation}
For Model {\sc i} a specific optimization program needs to be evaluated; it is relatively straightforward, but we leave out its specific form here.

\vb

{\it Organization. }
Section~\ref{NRARE} considers the situation in which the probability $p^{(N)}_t(a)$ does {\it not} correspond to a rare event (i.e., does not vanish as $N\to\infty$); the result is in terms of the distribution of the Poisson parameter of $M^{(N)}(t)$ (of which we characterize the density in terms of a system of partial differential equations).  {I}n Section \ref{AROUND} {we }study the distribution of $\phi_t(J)$ and $\psi_t(J)$ for values close to the maximum values they can attain (i.e., $a_t^\api$ and $a_t^\apii$). 
These results are then used in Section~\ref{RARE}, which covers the case in which $p^{(N)}_t(a)$  decays essentially exponentially as $N\to\infty$;
along the lines described above, we determine the exact asymptotics. Section \ref{numsec} contains remarks on computational aspects, as well as a set of numerical experiments. The paper is concluded by a discussion of the results obtained in Section \ref{CONCL}.

\section{Exact asymptotics in `non-rare range' --- distribution of the Poisson parameter}
This section studies the behavior of the Poisson parameters $\phi_t(J)$ and $\psi_t(J)$ in detail. In the first subsection 
the obtained results are used to evaluate the asymptotics of $p_t^{(N)}(a)$ for $N$ large for the case that $a$ is smaller than 
$a_t^\api$ (for Model {\sc i}) or $a_t^\apii$ (for Model {\sc ii}). The second subsection focuses on the shape of the distribution just below $a_t^\api$ (resp.\ $a_t^\apii$).

\subsection{Exact asymptotics in non-rare range}
\label{NRARE}
We start by considering the situation that the event of interest is not increasingly rare as $N\to \infty.$ 
For the moment we focus on Model {\sc i}, where it is noted that a similar line of reasoning, {\it mutatis mutandis}, applies to Model {\sc ii}. If  {$\phi_t(J)> a$}, then evidently the  probability
that ${\mathbb P}(P(N\phi_t(J)) \ge Na)$ converges to 1 as $N\to\infty$, and otherwise to 0. As a consequence,
\[\lim_{N\to\infty} p^{(N)}_t(a) = {\mathbb P}\left(\phi_t(J) \ge a\right).\]
As a consequence, we wish to  characterize the probabilities
${\mathbb P}(\phi_t(J) \ge a),$ and ${\mathbb P}(\psi_t(J) \ge a)$; the main result of this section is a system of partial differential equations that enables the evaluation of these objects.
For ease we assume that there are no distinct $i,j$ such that both $\lambda_i=\lambda_j$ and $\mu_i=\mu_j$; we comment later, in Remark \ref{dub}, on how to relax this assumption.

\subsubsection*{Model {\sc i}}
Our objective is to characterize the quantity \[p_i(a,t):= {\mathbb P}(\phi_t(J) \ge a,J(t)=i),\] for $i\in\{1,\ldots,d\}$, where it is assumed that $J(0)=i_0\in\{1,\ldots,d\}$. Consider the last $\Delta>0$ time units immediately before time $t$,  {$\Delta$} to be typically thought of as a small number. In this time interval the background process either jumps to state $i$ from a state $j\not=i$, or it was already in state $i$;  the third option, corresponding with  two or more jumps,  has probability $o(\Delta)$. 

If the process does not jump, then
\begin{eqnarray*}
\phi_t(J)&=& \int_0^{t-\Delta} \lambda_{J(s)} e^{-\int_s^t \mu_{J(r )}{\rm d}r} {\rm d}s +\:\int_{t-\Delta}^t\lambda_{\JGB{i}\ERASE{J(s)}} e^{-\mu_i(t-s)}{\rm d}s\\
&=&e^{-\mu_i\Delta } \int_0^{t-\Delta} \lambda_{J(s)} e^{-\int_s^{t-\Delta} \mu_{J(r )}{\rm d}r} {\rm d}s+\:\lambda_i \Delta+o(\Delta)\\
&=&(1-\mu_i\Delta ) \int_0^{t-\Delta} \lambda_{J(s)} e^{-\int_s^{t-\Delta} \mu_{J(r )}{\rm d}r} {\rm d}s+\:\lambda_i \Delta+o(\Delta),
\end{eqnarray*}
which is $(1-\mu_i\Delta ) \phi_{t-\Delta}(J)+\lambda_i \Delta+o(\Delta).$
As a consequence, up to terms of order $o(\Delta)$,
\begin{eqnarray*}
p_i(a,t)=\sum_{j\not =i} q_{ji}\Delta\, p_j(a,t) +\:\left(
1-\sum_{j\not=i} q_{ij}\Delta\right) p_i( a-\lambda_i\,\Delta+a\mu_i\,\Delta,t-\Delta).
\end{eqnarray*}
Subtracting $p_i(a,t)$ from both sides, dividing by $\Delta$, and letting $\Delta\downarrow 0$ leads to the following system of partial differential equations, for $i=1,\ldots,d$:
\[\sum_{j=1}^d q_{ji}p_j(a,t) =\frac{\partial}{\partial t}p_i(a,t) + (\lambda_i-a\mu_i) \frac{\partial}{\partial a}p_i(a,t).\]
We thus arrive at the following result; we present it in a compact form by using self-evident vector/matrix notation. 

\begin{prop} \label{p1} Consider Model {\sc i}. Assume $ a_t^\ami\le a\le a_t^\api.$
As $N\to\infty$,
\[ p^{(N)}_t(a) \to {\mathbb P}\left(\phi_t(J) \ge a\right)=\sum_{i=1}^d p_i(a,t) ,\]
where ${\bs p}(a,t)$ solves the system of partial differential equations
\[ Q^{\rm T} {\bs p}(a,t)= \frac{\partial}{\partial t}{\bs p}(a,t)+(\Lambda - a {\mathcal M})\frac{\partial}{\partial a}{\bs p}(a,t) .\]
\end{prop}

Now focus on additional conditions that are to be imposed. Recall that $J(0)=i_0$. \begin{itemize}
\item Let us start by identifying the conditions related to $t=0$. Realizing that $a_0^\ami=a_0^\api=0$, we have that $p_{i_0}(0,0) = 1$ and $p_i(0,0)=0$ for $i\not=i_0$. 
\item
Now consider the $a$-related conditions. Observe that 
\begin{eqnarray*}
{\mathbb P}\left(\phi_t(J) = \int_0^t \lambda_{i_0} e^{-\int_s^t \mu_{i_0}{\rm d}r}{\rm d}s\right)={\mathbb P}\left(\phi_t(J)= \frac{\lambda_{i_0}}{\mu_{i_0}} \left(1-e^{-\mu_{i_0} t}
\right)\right) = e^{-q_{i_0}t}.
\end{eqnarray*}
It follows that
\[p_i\left(a_t^\ami,t\right) = (e^{Qt})_{i_0,i},\:\:\:p_i\left(a_t^\api,t\right)=0\] for all $i\in\{1,\ldots,d\}$, but $p_{i_0}({\mathbf{\cdot}},t)$ has the special feature of having an atom of size $e^{-q_{i_0}t}$ at the value
\[a_t\s:=  \frac{\lambda_{i_0}}{\mu_{i_0}} \left(1-e^{-\mu_{i_0} t}
\right)\in\left[a_t^\ami,a_t^\api\right].\]
\end{itemize}

\begin{rem} \label{dub} Above we imposed the assumption that there are no distinct $i,j$ such that both $\lambda_i=\lambda_j$ and $\mu_i=\mu_j$. We now sketch what to do when this property does not hold. Let us consider the case that there is precisely one $j\not=i_0$ such that both $\lambda_{i_0}=\lambda_j$ and $\mu_{i_0}=\mu_j$; further generalizations can be performed in the same manner.  It is noted that now the atom at $a_t\s$ has size
\begin{eqnarray*}
e^{-q_{i_0}t} + \int_0^t  q_{i_0}e^{-q_{i_0}s} \cdot \frac{q_{i_0j}} {q_{i_0}} \cdot
e^{-q_j (t-s)}{\rm d}s=e^{-q_{i_0}t}+ \frac{e^{-q_jt}-e^{-q_{i_0}t}}{q_{i_0}-q_j} q_{i_0j}.\end{eqnarray*}
\end{rem}

\subsubsection*{Model {\sc ii}}
For Model {\sc ii} a similar approach can be followed.
We now concentrate on the object \[\bar p_i(a,t):= {\mathbb P}(\psi_t(J) \ge a\,|\,J(0)=i),\] for $i\in\{1,\ldots,d\}$.
Observe the subtle difference with the analysis of Model~{\sc i}: where we there considered the distribution of 
$\phi_t(J)$ {\it jointly with} $J(t)=i$, we now study the distribution of $\psi_t(J)$ {\it conditional on} $J(0)=i$.

Consider the first $\Delta>0$ time units, in which the background process either jumps, or stays in state $i$ (or jumps twice or more, but this corresponds to a probability that is $o(\Delta)$). 
If the process does not jump in $(0,\Delta]$, then, in distribution,
\begin{eqnarray*}
\psi_t(J) &=& \int_0^\Delta \lambda_i e^{- (t-s) {\mu_i}}{\rm d}s+\: \int_\Delta^t \lambda_{J(s)} e^{-(t-s) \mu_{J(s)}}{\rm d}s
\\& {\stackrel{\rm d}{=}}& \lambda_ie^{-\mu_i t}\,\Delta+\:\int_0^{t- \Delta}\lambda_{J(s)} e^{-(t-\Delta-s) \mu_{J(s)}}{\rm d}s+o(\Delta) ,\end{eqnarray*}
which is $ \lambda_ie^{-\mu_i t}\,\Delta+\psi_{t-\Delta}(J)+o(\Delta)$.
We thus find that
\[\bar p_i(a,t)= \sum_{j\not=i} q_{ij}\Delta\,\bar p_j(a,t)+\,\left(
1-\sum_{j\not=i} q_{ij}\Delta\right) \bar p_i( a-\lambda_ie^{-\mu_i t}\,\Delta,t-\Delta)+o(\Delta).\]
We continue in the usual way: subtracting $\bar p_i(a,t)$ from both sides, dividing by $\Delta$, and letting $\Delta\downarrow 0$ leads to the following system of partial differential equations, for $i=1,\ldots,d$:
\[\sum_{j=1}^d q_{ij}\bar p_j(a,t) =\frac{\partial}{\partial t}\bar p_i(a,t) + \lambda_ie^{-\mu_i t} \frac{\partial}{\partial a}\bar p_i(a,t).\]
This leads to the following statement, again in self-evident notation. 
\begin{prop} Consider Model {\sc ii}. Assume $a_t^\amii\le a\le a_t^\apii.$
As $N\to\infty$,
\[ p^{(N)}_t(a) \to {\mathbb P}\left(\psi_t(J) \ge a\right)=\sum_{i=1}^d \bar p_i(a,t) ,\]
where ${{\bar{\bs p}}}(a,t)$ solves the system of partial differential equations
\[Q\,\bar{\bs p}(a,t) = \frac{\partial}{\partial t}\,\bar{\bs p}(a,t)+(\Lambda\, e^{-{\mathcal M}t})\frac{\partial}{\partial a}\,\bar{\bs p}(a,t).\]
\end{prop}

Again additional conditions should be imposed:
\begin{itemize}
\item We have $\bar p_i(0,0)=1$ for all  $i\in\{1,\ldots,d\}.$
\item In this case $\bar p_i(a_t^-,t) = 0$ and $\bar p_i(a_t^+,t)=1$  {for all $i\in\{1,\ldots,d\}$ and}  $p_{i}(a,t)$ has  an atom of size $e^{-q_{i}t}$ at 
\[a_{i,t}\s:=  \frac{\lambda_{i}}{\mu_{i}} \left(1-e^{-\mu_{i} t}
\right).\]\end{itemize}
It is noted that these conditions can be adapted in case there is a $j$ such that $\lambda_{i_0}=\lambda_j$ and $\mu_{i_0}=\mu_j$, in the way pointed out in Remark \ref{dub}.

\subsection{Distribution of Poisson parameter close to its domain {boundaries}}\label{AROUND}
In this section we study the behavior, for small $\delta${,} of $\phi_t(J)$ and $\psi_t(J)$ being less than $\delta$ away from $a_t^\api$ and $a_t^\apii$, respectively.
The exposition is slightly easier for Model {\sc ii}, due to the fact that for that model the maximum attainable variable is explicitly known (see (\ref{MAXP})), but for Model {\sc i} essentially the same approach can be followed. The results obtained in this subsection are crucial when deriving the exact asymptotics in Section \ref{RARE}.

Define the `maximizing path'
\[\gamma_t(s):=\arg\max_{i\in\{1,\ldots,d\}} \lambda_i e^{-(t-s)\mu_i}.\]
As was shown in \cite{BKMT,BM} $\gamma_t(\cdot)$ jumps at most $d-1$ times in $[0,t]$; let $D\leq d-1$ be this number of jumps. Then there are two cases: no jumps at all in $[0,t]$, and a positive number of jumps {in $[0,t]$ }(in which case we denote by $s_1$ up to $s_D$ the epochs of these jumps). The former case being elementary,  we focus in this section on the latter case. Without loss of generality we assume that the states are labeled such that $\gamma_t(s)$ visits the states $1$ up to $D+1$ when $s$ increases from $0$ to $t$.

We first evaluate the difference between the maximum value $a_t^\apii$ of $\psi_t(J)$ (corresponding to jumps at $s_1$ up to $s_D$) with the value of $\psi_t(J)$
that results from jumps at times $s_1+v_1\varepsilon$ up to $s_D+v_D\varepsilon${, where the $v_i\varepsilon$ are small (but not necessarily positive)}. It is readily checked that this difference equals, with $s_0=0$, $s_{D+1}=t$, and $v_0=v_{D+1}=0$,
\[\sum_{i=1}^{D+1}\left(\int_{s_{i-1}}^{s_i} \lambda_i e^{-\mu_i(t-r)}{\rm d}r -\int_{s_{i-1}+v_{i-1}\varepsilon}^{s_i+v_{i}\varepsilon} \lambda_i e^{-\mu_i(t-r)}{\rm d}r\right),\]
which can alternatively be written as
\[\sum_{i=1}^{D+1} 
\frac{\lambda_i}{\mu_i} e^{-\mu_i t} \left(e^{\mu_i s_i} - e^{\mu_i s_{i-1}}\right)-\sum_{i=1}^{D+1} 
\frac{\lambda_i}{\mu_i} e^{-\mu_i t} \left(e^{\mu_i (s_i+v_i\varepsilon)} - e^{\mu_i (s_{i-1}+v_{i-1}\varepsilon)}\right),
\]
or, further simplified,
\begin{equation}\sum_{i=1}^{D+1}\frac{\lambda_i}{\mu_i} e^{-\mu_i(t-s_i)} (1-e^{\mu_i v_i\varepsilon})
-\sum_{i=1}^{D+1}\frac{\lambda_i}{\mu_i} e^{-\mu_i(t-s_{i-1})} (1-e^{\mu_i v_{i-1}\varepsilon}),\label{func}\end{equation}
notice that, due to $v_0=v_{D+1}=0$ the last term of the first sum can be left out, and the same holds for the first term of the second sum.
Recalling that, immediately from the definition of $s_1,\ldots,s_D$,
\[\lambda_i e^{-\mu_i(t-s_i)} = \lambda_{i+1} e^{-\mu_{i+1}(t-s_i)},\:\:i=1,\ldots,D,\]
we have that (\ref{func}) equals,   up to terms that are $o(\varepsilon^2)$, 
\[\sum_{i=1}^D\left(\lambda_{i}  e^{-\mu_{i }(t-s_i)}- \lambda_{i{+1}}  e^{-\mu_{i{+1}}(t-s_i)}  \right)v_i\varepsilon
+\sum_{i=1}^D\omega_i\, (v_i\varepsilon)^2= \sum_{i=1}^D\omega_i\, (v_i\varepsilon)^2;\:\:
\]
here we have used the definition, for $i=1,\ldots,D$,
\begin{eqnarray*}\omega_i&:=&{   {\frac{\lambda_{i+1}\mu_{i+1}}{2}} e^{-\mu_{i+1}(t-s_i)} -
 {\frac{\lambda_{i}\mu_{i}}{2}} e^{-\mu_{i}(t-s_i)}}\\&=&\frac{\lambda_{i+1}}{2}(\mu_{i+1}-\mu_i) e^{-\mu_{i+1}(t-s_i)}=\frac{\lambda_{i}}{2}(\mu_{i+1}-\mu_i) e^{-\mu_{i}(t-s_i)}.\end{eqnarray*}
 
  It is readily verified that along $\gamma_t(\cdot)$ it holds that $\mu_i\geq\mu_j$ if $i>j$, and hence all coefficients $\omega_i$ are non-negative; this is in line with the fact that the functional $\psi_t(J)$ is maximized  by the path $\gamma_t(\cdot)$.  
We thus arrive at
\[{\mathbb P}\left(\psi_t(J) \ge a^\apii_t - {\delta}\right)= \pi_1 q_{1} e^{-q_1 s_1} \frac{q_{12}}{q_1}
{q_2}e^{-q_2 (s_2-s_1)} 
\cdots \frac{q_{D,D+1} }{q_D} {q_{D+1}}e^{-q_{D+1} (t-s_D)} {\mathscr V}(\delta)+o({\mathscr V}(\delta))
,\]
where ${\mathscr V}(\delta)$ denotes the volume of the set
\[{\mathscr S}(\delta):=\left\{ (x_1,\ldots, x_D):\,\sum_{i=1}^D\omega_i\,x_i^2 <\delta\right\},\] which is $\kappa_t\cdot R^D=\kappa_t\cdot{\delta}^{D{/2}}$ for some constant $\kappa_t>0$ and $R:=\sqrt{\delta}$ being the `scale' of the ellipsoid.
We have thus identified a constant $\bar \kappa_t>0$ such that
\[\lim_{\delta\downarrow 0} {\mathbb P}\left(\psi_t(J) \ge a^\apii_t - {\delta}\right) {\delta}^{-D {/2}} =\bar \kappa_t.\]
A similar argument provides us with the corresponding density close to $a^\apii_t$; then essentially the integration needs to be done over $\partial {\mathscr S}(\delta)$, which is of the order 
$R^{D-1}$. Appealing to the chain rule (with ${\rm d}R/{\rm d}\delta = (2\sqrt{\delta})^{-1}$), we thus find that for a constant $\hat \kappa_t>0$,
\begin{equation}
\label{esti}\lim_{\delta\downarrow 0} {\mathbb P}\left(a^\apii_t-\psi_t(J) \in {\rm d}{\delta}\right) \delta^{-D/2+1} =\hat \kappa_t.\end{equation}

We note that above we tacitly   imposed the regularity condition that all transition rates along the path $\gamma_t(\cdot)$ are positive:
\begin{equation}
\label{reg} q_{i,i+1}>0\:\:\mbox{for all $i\in\{1,\ldots,D\}$.}
\end{equation}
As an aside we mention that adaptation of the arguments to the case in which along $\gamma_t(\cdot)$ there are (one or more) states $i\in\{1,\ldots,D\}$   corresponding with $q_{i,i+1}=0$ is a purely technical issue, and is relatively straightforward. Importantly, it can be checked that it affects the power of $\delta$ appearing in (\ref{esti}).  {E}xample 2 illustrate{s} how this issue can be dealt with. 

\begin{exmp} \label{EX1} Consider Model {\sc ii} with $d=2$. We consider the case that $\lambda_1<\lambda_2$ and $\mu_1<\mu_2$, so that the curves $\lambda_i e^{-\mu_i(t-s)}$ intersect at
\[s{_1}= t -\bar s,\:\:\:\mbox{with}\:\:\:\bar s:=\frac{\log(\lambda_1/\lambda_2)}{\mu_1-\mu_2};\]
we assume $t>\bar s$. Because of the choice of our parameters, we are in the situation that the maximizing path jumps once {in $[0,t]$},
where
\[\omega_1 = \frac{\lambda_2}{2} (\mu_2-\mu_1) e^{-\mu_2(t-s{_1})}= \frac{\lambda_2}{2}(\mu_2-\mu_1) \left(\frac{\lambda_1}{\lambda_2}\right)^{-\mu_2/(\mu_1-\mu_2)}.\]
We conclude that
\[{\mathscr V}(\delta) = \frac{ {2\sqrt{2\delta}}}{\sqrt{\lambda_2(\mu_2-\mu_1)}} \sqrt{ \left(\frac{\lambda_1}{\lambda_2}\right)^{\mu_2/(\mu_1-\mu_2)}},\]
and hence 
\[{\bar\kappa_t} = \pi_1\,{q_{12}q_2}\,e^{-q_1 t} \frac{ 2\sqrt{2}}{\sqrt{\lambda_2(\mu_2-\mu_1)}}    \left(\frac{\lambda_1}{\lambda_2}\right)^{(q_1-q_2+\mu_2/2)/(\mu_1-\mu_2)}.\]
\end{exmp}

\begin{exmp} 
In this  example we consider a situation in which  regularity condition  (\ref{reg})  does not apply. We point out how in this case the density close to $a_t^\apii$ can be evaluated. As becomes clear, the procedure is straightforward but tedious; therefore we assume in the next section, when evaluating the asymptotics, that the simpler situation in which (\ref{reg}) is in place.

We consider the same setting as in the previous example, but now with $d=3$ where the transition rates $q_{ij}$ are such that state 2 can be reached from state 1 {\it only via state 3}: $q_{13}, q_{32}>0$ but $q_{12}= 0$. We assume that for any $s\in[0,t]$ the function $\lambda_3e^{-\mu_3(t-s)}$ nowhere majorizes $\lambda_1e^{-\mu_1(t-s)}$ or $\lambda_2e^{-\mu_2(t-s)}$. In other words: as in the previous example the maximizing path subsequently visits states 1 and 2 (and the resulting value of $a_t^\apii$ is the same), but the  modulating Markov chain cannot jump directly from state 1 to 2.

Consider the path at which there is a transition from state 1 to 3 at time $s{_1}-v_1\varepsilon$, and then a transition from state 3 to 2 at time $s{_1}+ v_2\varepsilon${, with $v_i\varepsilon$ small and positive}. The difference between $a_t^\apii$ and the value of $\psi_t(J)$ resulting from this path is
\[\frac{\lambda_1}{\mu_1}e^{-\mu_1(t-s{_1})} \left(1-e^{-\mu_1v_1\varepsilon}\right)
+\frac{\lambda_2}{\mu_2}e^{-\mu_2(t-s{_1})} \left(e^{\mu_2v_2\varepsilon}-1\right)-
\frac{\lambda_3}{\mu_3}e^{-\mu_3(t-s{_1})} \left(e^{\mu_3v_2\varepsilon}-e^{{-}\mu_3v_1\varepsilon} \right),\]
which behaves, for $v_i\varepsilon$ small, as
$z_1 v_1\varepsilon +z_2 v_2\varepsilon$, with
$z_i:= \lambda_ie^{-\mu_i(t-s{_1})} - \lambda_3e^{-\mu_3(t-s{_1})};$ recall that $z_i>0.$
We thus arrive at, ignoring terms that are $o({\mathscr V}(\delta))$,
\begin{eqnarray*}{\mathbb P}\left(\psi_t(J) \ge a^\apii_t - {\delta}\right)&=& 
\pi_1 q_{1} e^{-q_1 s {_1} }\frac{q_{13}}{q_1} q_3 e^{-q_3\cdot 0}
\frac{q_{32}}{q_3} {q_2}e^{-q_{{2}} (t-s{_1})} {\mathscr V}(\delta)\\&=&
\pi_1 q_{1{3}} e^{-q_1 s{_1} } {q_{32}}{q_2}e^{-q_{2} (t-s{_1})} {\mathscr V}(\delta)
,\end{eqnarray*}
where ${\mathscr V}(\delta)$ denotes the volume of the set 
\[{\mathscr S}(\delta):=\left\{ (x_1,x_2)\in{\mathbb R}_+^2:\,z_1\,x_1+ z_2\,x_2 <\delta\right\},\] i.e., $\delta^2/({2}z_1z_2)$.
Conclude that for $\delta$ small the probability under investigation is essentially proportional to $\delta^2$. This is in contrast with the order ${\sqrt{\delta}}$ that we found in Example 1; apparently the likelihood of reaching values close to $a^\apii_t$ is considerably smaller in Example 2, as a consequence of the additional transitions needed.
\end{exmp}

\section{Exact asymptotics in `rare range'}\label{RARE}
In the previous section we have considered the situation in which $p^{(N)}_t(a)$ converges to a positive constant; this case corresponds to the exceedance level $a$ being between the minimum and maximum value of the Poisson parameter underlying the distribution of $M^{(N)}(t)$. In the present section we look at the opposite case, i.e., the case in which $p^{(N)}_t(a)$ vanishes as $N$ grows large. We present the analysis for Model {\sc i}, but Model {\sc ii}  can be dealt with fully analogously. 

Below we consider the situation that $a>a^\api$; the asymptotic analysis of $1- p^{(N)}_t(a)$ for $a< a^\ami$ follows in the same way.
To this end, we first realize that we have the following representation, due to the fact that $M^{(N)}(t)$ has a Poisson distribution with random mean:
\[p^{(N)}_t(a)=  \int_{a_t^\ami}^{a_t^\api} {\mathbb P}(P(N\alpha)\ge Na)\,
{\mathbb P}(\phi_t(J)\in {\rm d}\alpha);\]
the integral is on the interval $[{a_t^\ami},{a_t^\api}]$, as this is the interval of values that
$\phi_t(J)$ can attain.

The first step is to analyze ${\mathbb P}(P(N\alpha)\ge Na)$, relying on standard probabilistic tools. 
Define, for $\alpha\in[{a_t^\ami},{a_t^\api}],$ with $\Lambda(\vt\,|\,\alpha):= \log {\mathbb E}\,e^{\vt P(\alpha)}$,
the {\it Legendre transform}
\begin{eqnarray*}
I(a\,|\,\alpha)&:=&\sup_{\vt}\left(\vt a -\Lambda(\vt\,|\,\alpha)\right)=\sup_{\vt}\left(\vt a -\alpha(e^\vt-1)\right).\end{eqnarray*}

As the optimizing $\vt$ equals $\vt(a\,|\,\alpha)= \log(a/\alpha)>0,$
we have
$I(a\,|\,\alpha) = a\log({a}/{\alpha})+\alpha-a.$
As can be found in e.g.\ \cite{DZ}, the lattice version of the Bahadur-Rao result \cite{BR} states that, as $N\to\infty$,
\[{\mathbb P}(P(N\alpha)\ge Na) \cdot \left(e^{N I(a\,|\,\alpha)} \sqrt{{2\pi}N} \cdot 
\xi(a\,|\,\alpha)\right)\to 1,\]
where
\begin{eqnarray*}
\xi(a\,|\,\alpha)&:=& \sqrt{\Lambda''(a\,|\,\alpha)}\left({1-e^{-\vt(a\, {|}\,\alpha)}}\right)=
\sqrt{a}\left(1-\frac{\alpha}{a}\right).
\end{eqnarray*}
Interestingly, we know that this convergence is {\it uniform} in $\alpha\in[{a_t^-},{a_t^+}],$ as an immediate consequence of the results in H\"oglund \cite{HOG}. This implies that, for all $\varepsilon>0$ we have that for $N$ large enough
\begin{eqnarray*}{\sup_{\alpha\in[{a_t^-},{a_t^+}]} 
{\mathbb P}(P(N\alpha)\ge Na) \cdot\left( e^{N I(a\,|\,\alpha)} \sqrt{{2\pi}N} \cdot 
\xi(a\,|\,\alpha)\right)}\in (1-\varepsilon, 1+ \varepsilon).\end{eqnarray*}
In addition, we have, uniformly in $N$, the celebrated {\it Chernoff bound}:
\begin{equation}
\label{ch}{\mathbb P}(P(N\alpha)\ge Na)\le e^{-N I(a\,|\,\alpha)}.\end{equation}
When analyzing the asymptotics of $p^{(N)}_t(a)$ for $N$ large and $a>a_t^\api$, two cases need to be distinguished: the case that $\phi_t(J)$ does not have an atom in $a_t^\api$, and the case that it has. Let us start with the former case (which is more involved than the latter case). 

\vb

$\rhd$ {\it Case $1$ --- $\phi_t(J)$ does not have an atom in $a_t^\api$.}
Fix some $\delta\in(-1,-\frac{1}{2})$.  We split  $p^{(N)}_t(a)$ into
\begin{equation}
\label{split}
K\left(a_t^\ami,a_t^\api-N^\delta\right)+
K\left(a_t^\api-N^\delta,a_t^\api\right),\end{equation}
where,
for $u<v$, \[K(u,v):=\int_{u}^{v} {\mathbb P}(P(N\alpha)\ge Na)\,
{\mathbb P}(\phi_t(J)\in {\rm d}\alpha).\]

Let us start by analyzing the first term in (\ref{split}); our goal is to show that it can be ignored (asymptotically, i.e., as $N\to\infty$) relative to the second term.
Observe that, because of (\ref{ch}),
\begin{equation}{e^{N I(a\,|\,a_t^\api)}\,K\left(a_t^\ami,a_t^\api-N^\delta\right)} \le  A_t\,
e^{N I(a\,|\,a_t^\api)} \left( \sup_{\alpha\in [a_t^\ami,a_t^\api-N^\delta]}  e^{-N I(a\,|\,\alpha)} \right)\,\label{T1}\end{equation}
where $A_t:= a_t^\api-a_t^\ami$; in view of the shape of the asymptotic expansion that eventually  comes out, we multiplied by   $e^{N I(a\,|\,a_t^\api)}$.
Now realize that $I(a\,|\,\alpha)$ is convex in $\alpha$, having the value $0$ when $\alpha=a$, and that it is decreasing in $\alpha$, since $a>a_t^\api$. It thus follows that
\[\arg\inf_{\alpha\in [a_t^\ami,a_t^\api-N^\delta]}I(a\,|\,\alpha)=a_t^\api-N^\delta.\]
As a consequence, (\ref{T1}) is majorized by
\begin{equation}\label{T1-b}A_t\,e^{N I(a\,|\,a_t^\api)}\,e^{-N I(a\,|\,\alpha_t^\api-N^\delta)}.\end{equation}
We now present an upper bound on the exponent featuring in (\ref{T1-b}). It is a trivial exercise to verify that standard estimates yield
\begin{eqnarray*}{I(a\,|\,a_t^\api) - I(a\,|\,{a}_t^\api-N^\delta)} =  a \log\frac{a_t^\api-N^\delta}{a_t^\api}+N^\delta
\le \left(1-\frac{a}{a_t^\api}\right)N^\delta \le -c N^\delta,\end{eqnarray*}
for some positive $c$ (where it is used that $a>a_t^\api$). Conclude that Expression (\ref{T1-b}) is bounded from above by
${A_t}\, \exp({-cN^{1-\delta}}),$ and therefore we obtain, as $N$ grows large,
\begin{eqnarray}{N^{(D+1)/2}\,e^{N I(a\,|\,a_t^\api)}\,K\left(a_t^\ami,a_t^\api-N^\delta\right)}
\le
N^{(D+1)/2} \,A_t\,e^{-cN^{1-\delta}}\label{Int1}
\to 0.\end{eqnarray}

Let us now concentrate on the second term in (\ref{split}); as we will show, it dominates the contribution of the first term. To this end, we first focus on an upper bound, but, as we see later on, a corresponding lower bound can be derived very similarly, thus establishing the exact asymptotics of $p_t^{(N)}(a)$. 
Because of the (uniform version of) the Bahadur-Rao result (as was stated above), we have that for any $\epsilon>0$,
\begin{eqnarray}\nonumber\lefteqn{\limsup_{N\to\infty}N^{(D+1)/2}e^{N I(a\,|\,a_t^\api)}\,K\left(a_t^\api-N^\delta,a_t^\api\right)}\\
&\le &(1+\varepsilon)\, {\cdot}\,\limsup_{N\to\infty}N^{D/2} \int_{a_t^\api-N^\delta}^{a_t^\api} G_N(\alpha)\,{\mathbb P}(\phi_t(J)\in{\rm d}\alpha),\label{rhs}
\end{eqnarray}
where, with $\eta(a\,|\,\alpha):=1/(\sqrt{2\pi}\;\xi(a\,|\,\alpha))$,
\[G_N(\alpha):= e^{N\,I(a\,|\,a^\api_t)-N\,I(a\,|\,\alpha)} \,\eta(a\,|\,\alpha).\]

We now further analyze (\ref{rhs}). To this end, we first define
\[\bar G(\alpha):=a\log\left(1-\frac{\alpha}{a_t^\api}\right)+\alpha,\]
and assume that the regularity condition  (\ref{reg}) applies. 
By virtue of standard continuity arguments it follows that in combination with (\ref{esti}), for all $\varepsilon'>0$, Expression (\ref{rhs}) is majorized by 
\begin{eqnarray*}\lefteqn{\hspace{-1.4cm}(1+\varepsilon')\, \eta(a\,|\,a_t^\api)\,\hat\kappa_t{\cdot}\,\limsup_{N\to\infty}N ^{D/2}
\int_{a_t^\api-N^\delta}^{a_t^\api} e^{N\,I(a\,|\,a^\api_t)-N\,I(a\,|\,\alpha)}(a^\api_t-\alpha)^{D/2-1}\,{\rm d}\alpha}\\
&\stackrel{{\beta\,:=\,a^\api_t-\alpha}}{=}&(1+\varepsilon')\, \eta(a\,|\,a_t^\api)\,\hat\kappa_t {\cdot}\,\limsup_{N\to\infty}N ^{D/2}
\int_{0}^{N^\delta}e^{N\bar G({\beta})}{\beta}^{D/2-1}{\rm d}{\beta} {.}\end{eqnarray*}
Using elementary Taylor expansions, it is easily verified that there are numbers $\ell$ and $u$ such that, with \[b:=\left(\frac{a}{a_t^\api}-1\right)>0,\] for $N$ sufficiently large and all ${\beta}\in[0,N^\delta]$,
\[\ell N^{1+2\delta} - b {\beta} N \le N\left(a\log\left(1-\frac{{\beta}}{a_t^\api}\right)+ {\beta}\right)
 \le u N^{1+2\delta} - b {\beta} N.\]
As a consequence,  using in step (i) that $\delta<-\frac{1}{2}$ and in step (ii) $\delta>-1$, 
\begin{eqnarray*}\limsup_{N\to\infty}N ^{D/2}
\int_{0}^{N^\delta} e^{N\bar G( {\beta})} {\beta}^{ {D/2-1}}{\rm d} {\beta}&\le& \limsup_{N\to\infty}N^{D/2}  e^{uN^{1+2\delta}} \int_0^{N^\delta} e^{-b {\beta} N} {\beta}^{D/2-1}{\rm d} {\beta} \\
&\stackrel{\rm (i)}{=}&
 \limsup_{N\to\infty}N^{D/2}  \int_0^{N^\delta} e^{-b {\beta} N} {\beta}^{D/2-1}{\rm d} {\beta} \\
 &\stackrel{ {\alpha\,:=\,b\beta N}}{=}&   \frac{1}{b^{D/2-1}}\limsup_{N\to\infty} \int_0^{bN^{\delta+1}} e^{- {\alpha}}  {\alpha}^{D/2-1} {\rm d} {\alpha}
\stackrel{\rm (ii)}{=}   \frac{\Gamma(D/2)}{b^{D/2}}.
\end{eqnarray*}
The corresponding lower bound can be found along the same lines: for an arbitrary $\varepsilon'>0$,
\begin{eqnarray*}\nonumber\lefteqn{\hspace{-1mm}\hspace{-1mm}\liminf_{N\to\infty}N^{(D+1)/2} e^{N I(a\,|\,a_t^\api)}\,K\left(a_t^\api-N^\delta,a_t^\api\right)}\\
&\ge& (1-\varepsilon')\, \eta(a\,|\,a_t^\api)\,\hat\kappa_t {\cdot}\, {\liminf}_{N\to\infty}N^{D/2}  e^{\ell N^{1+2\delta}} \int_0^{N^\delta} e^{-b\alpha N}\alpha^{D/2-1}{\rm d}\alpha,
\end{eqnarray*}
which can be evaluated as before.
By taking $\varepsilon'\downarrow 0$, upon combining the above upper and lower bound, we obtain 
\begin{eqnarray}{ \lim_{N\to\infty}N^{(D+1)/2}\,e^{N I(a\,|\,a_t^\api)}\,K\left(a_t^\api-N^\delta,a_t^\api\right)}= \eta(a\,|\,a_t^\api)\,\hat\kappa_t \frac{\Gamma(D/2)}{b^{D/2}}.\label{Int2} 
\end{eqnarray}
Next we combine the asymptotics of both intervals, i.e., the one over $[a_t^\ami,a_t^\api-N^\delta)$ and the one over $[a_t^\api-N^\delta, a_t^\api]$. 
From (\ref{Int1}) and (\ref{Int2}), the main result of this section follows. The analogous result for Model {\sc ii} can be derived in precisely the same way; the only difference lies in the value of the constant $\hat\kappa_t$.

\begin{thm} \label{st} Consider Model {\sc i}. Assume $a> a_t^\api$, and let $\phi_t(J)$ have no atom in $a_t^\api$; in addition, assume that regularity condition $(\ref{reg})$ applies. As $N\to\infty$, 
\[N^{(D+1)/2} \,e^{N I(a\,|\,a_t^\api)}\,p_t^{(N)}(a)\to \left(\frac{a_t^\api}{a-a_t^\api} \right)^{D/2}\frac{\hat\kappa_t\,\Gamma(D/2)}{ {\sqrt{2\pi}}\,\xi(a\,|\,a_t^\api)}.\]
\end{thm}

{\it $\rhd$ Case $2$ --- $\phi_t(J)$ has an atom in $a_t^\api$}. We now consider the situation that \[F(a_t^\api):={\mathbb P}\left(\phi_t(J) = a_t^\api\right)> 0.\] Because of the arguments used in the derivation of Thm.\ \ref{st}, we observe that the contribution to the probability of interest  due to the event $\phi_t(J)\in[a_t^\ami,a_t^\api)$ is of an order of at most \[\frac {e^{-N\, I(a\,|\,a_t^\api)}}{N}\]
(up to a multiplicative constant); 
realize that this is a consequence of the fact that the corresponding path requires at least one jump. From the Bahadur-Rao result, however, it is directly seen that the contribution due to the event $\phi_t(J)=a_t^\api$ is larger, viz.\ of the order (up to a multiplicative constant)  \[\frac{e^{-N\, I(a\,|\,a_t^\api)}}{\sqrt{N}}.\] As a consequence, the latter scenario dominates, and we obtain the following exact asymptotics; again, an analogous result is valid for Model {\sc ii}.

\begin{cor}\label{cor1}
Consider Model {\sc i}. Assume $a> a_t^\api$, and let $\phi_t(J)$ have an atom in $a_t^\api$. As $N\to\infty$, 
\[\sqrt{N} \,e^{N I(a\,|\,a_t^\api)}\,p_t^{(N)}(a)
\to\frac{F(a_t^\api)}{ {\sqrt{2\pi}}\,\xi(a\,|\,a_t^\api)}.\]\end{cor}

\section{Computational issues} \label{numsec}
The objective of this section is to present an efficient simulation method for estimating $p_t^{(N)}(a)$ for the situation that $a$ is large than (in Model {\sc i}) $a_t^\api$ or (in Model {\sc ii}) $a_t^\apii$. In addition we include a numerical experiment featuring a typical example.

\vb

{\it Basic method, and its logarithmic efficiency.} Particularly when $N$ is large, the probability $p_t^{(N)}(a)$ will be small, thus imposing constraints on the feasibility of standard Monte Carlo techniques. There is, however, an interesting remedy. To this end, note that we can express the probability of our interest as
\begin{equation}
\label{form}p_t^{(N)}(a) = {\mathbb E} {\mathscr P} (Na, N\phi_t(J))\end{equation}
 (where, as an aside, we mention that we point the procedure out for Model {\sc i}, but Model {\sc ii} can be dealt with fully analogously);
the function
\[{\mathscr P} (n, \lambda):=
\sum_{k=n}^\infty e^{-\lambda}\frac{\lambda^k}{k!},\]
is the tail distribution of the Poisson distribution, and is available in standard software packages.
The form (\ref{form})  suggests the following  simple and effective simulation approach: in  run ${{\ell}}$ (with ${{\ell}}=1,\ldots, {M}$) the path $J_{{\ell}}$ is sampled, the parameter $\phi_t(J_{{\ell}})$ is calculated, and the probability $p_t^{(N)}(a)$ is estimated by
\[\frac{1}{ {M}}\sum_{{{\ell}}=1}^{ {M}}  {\mathscr P} (Na, N\phi_t(J_{{\ell}})).\]
This procedure is {\it logarithmically efficient} \cite[Ch. VI]{ASMGL}. To see this, first note that we have the obvious deterministic upper bound
\begin{equation}\label{UB}{\mathscr P} (Na, N\phi_t(J))\le  {\mathscr P} (Na, N a_t^\api),\end{equation} 
as a consequence of the stochastic monotonicity of the Poisson distribution in its parameter.
Due to  Jensen's inequality in combination with Thm.\ \ref{st} and Corollary \ref{cor1} we have the lower bound
\[\liminf_{N\to\infty}\frac{1}{N}\log {\mathbb E}  {\mathscr P}^2 (Na, N\phi_t(J)) \ge
2\,\lim_{N\to\infty}\frac{1}{N}\log {\mathbb E}  {\mathscr P} (Na, N\phi_t(J)) = -2 I(a\,|\,a_t^\api).\]
Because of (\ref{UB}), however, this lower bound is actually achieved:
\[\limsup_{N\to\infty}\frac{1}{N}\log {\mathbb E}  {\mathscr P}^2 (Na, N\phi_t(J)) \le
2\,\lim_{N\to\infty}\frac{1}{N}\log  {\mathscr P} (Na, N a_t^\api) = -2 I(a\,|\,a_t^\api).\]
We thus obtain logarithmic efficiency. Often simulation experiments are performed until the estimate has reached a certain efficiency: the ratio of the width of the confidence interval to the estimate is smaller than some predefined number (e.g.\ 10\%). In practical terms, in this setting with $p_t^{(N)}(a)$ decaying essentially exponentially in $N$, logarithmic efficiency effectively means that the number of runs that is needed grows at most {\it subexponentially} in $N$.

\vb

{\it Importance-sampling based acceleration.}
In fact, the rare event studied in this paper is the effect of the combination of (i) the Poisson parameter $\phi_t(J)$ attaining a rare value, say $\phi$, and (ii) a Poisson random variable with parameter $N\phi$ attaining a rare value. Note that the above approach adequately deals with the randomness due to effect (ii) -- that is, we do not need to sample the Poisson random variable, but we use computations instead.

The question that is left concerns the rarity which is a consequence of $\phi_t(J)$ attaining a rare value. In the proofs we have seen that overflow is most likely caused by $\phi_t(J)$ attaining a value `close to' its maximal value $a_t^\api$, which only happens when the jump epochs are close to those of some maximizing path (that was explicitly determined in \cite{BKMT} and \cite{BM} for Models {\sc i} and {\sc ii}, respectively). We saw that the probability of $\phi_t(J)$ being an amount in the order of $\delta$ away from its maximum value $a_t^\api$, is of the order $\delta^{D/2}$, i.e., relatively rare. Importance sampling can be used to resolve this issue in the following way.

Choose $\Delta$ sufficiently {small} such that all $s_i$ pairs are at least $2 \Delta$ apart; recall that the $s_i$ are the transition epochs along the path that optimises the Poisson parameter. We let $T_0=0$ and $T_i$, for $i=1,2,\ldots,  {D+1}$ be the subsequent transition epochs of the background process in our simulation, and $U_i:=T_i-T_{i-1}$ the corresponding sojourn times.
We write, with $\gamma(\cdot)$ being functions that map $[0,t]$ onto $\{1,\ldots,d\}$ and $\bar s_i:=s_i-s_{i-1}$,
\[{\mathscr Z}(\Delta) := \left\{
\gamma(\cdot)\left| 
\begin{array}{l}\gamma(s) = i \:\:\mbox{$\forall s\in [T_{ {i-1}},T_{i})\:{\:\:\forall i=1,\ldots,D+1}$};
\\ U_i \in (\bar s_i-\Delta,\bar s_i+\Delta)\:\:\mbox{$\forall i=1,\ldots,D$};
\\ 
{U_{D+1}\ge t - s_D + D\Delta}
\end{array}\right.
\right\}.\]
The set ${\mathscr Z}(\Delta) $ should be interpreted as the collection of paths that are `close to' the path that maximizes the random parameter of the Poisson distribution; recall that, without loss of generality, we had labeled the states such that along this optimizing path the states $1$ up to $D+1$ are subsequently visited.

The idea is now to estimate the quantities 
\[ {\mathbb E} \left({\mathscr P} (Na, N\phi_t(J)) \,1\{J\not\in {\mathscr Z}(\Delta)\}\right)\:\:\:\mbox{and}\:\:\:
{\mathbb E} \left({\mathscr P} (Na, N\phi_t(J)) \,1\{J\in {\mathscr Z}(\Delta)\}\right)\]
separately, and to add the resulting estimates up. The first of these quantities is estimated
under the actual measure ${\mathbb P}$, whereas for the second (which contains the rare event 
of $\phi_t(J)$ being close to $a_t^\api$) we use importance sampling. In more detail:
\begin{itemize}
\item[$\circ$]
The quantity ${\mathbb E} \left({\mathscr P} (Na, N\phi_t(J)) \,1\{J\not\in {\mathscr Z}(\Delta)\}\right)$ is estimated
by performing ${M}_1$ runs: 
\[\frac{1}{{M}_1}\sum_{{{\ell}}=1}^{{M}_1}  {\mathscr P} (Na, N\phi_t(J_{{\ell}})) \,1\{J_i\not\in {\mathscr Z}(\Delta)\},\]
with the $J_{{\ell}}$ sampled under ${\mathbb P}.$
\item[$\circ$]
The quantity ${\mathbb E} \left({\mathscr P} (Na, N\phi_t(J)) \,1\{J\in {\mathscr Z}(\Delta)\}\right)$ can be estimated using an
importance sampling approach: an alternative measure, say ${\mathbb Q}$, is used to draw samples $\phi_t(J_1)$ up to $\phi_t(J_{{M}_2})$, and then the simulation output (i.e., ${\mathscr P} (Na, N\phi_t(J_{{\ell}}))$) is translated back in terms of the original probability measure ${\mathbb P}$ by multiplying it with an appropriate likelihood ratio $L_{{\ell}}$ (to be interpreted as a Radon-Nikodym derivative ${\rm d}{\mathbb P}/{\rm d}{\mathbb Q}$).

The measure ${\mathbb Q}$ is constructed as follows. The transition probabilities are changed in such a way that with probability 1 the background process visits the states $1$ up to $D+1$. Along this path, the time spent in state $i$ is sampled from a distribution with density, for $s\in (\bar s_i -\Delta,\bar s_i + \Delta)$,
\[ {q_ie^{-q_i s} } \left(\int_{\bar s_i-\Delta}^{\bar s_i+\Delta} q_ie^{-q_i r}{\rm d}r\right)^{-1}= \frac{q_ie^{-q_i s}}{
\sigma_i},\:\:\:\mbox{with}\:\:\:\sigma_i:=e^{-q_i(\bar s_i-\Delta)}-e^{-q_i(\bar s_i+\Delta)} \]
(where the density is defined to be $0$ elsewhere), for $i=1,\ldots,D$. The time spent in state $D+1$ is sampled from a distribution with density, for $s\ge  t-s_D+ D\Delta $,
\[ \frac{q_{D+1}e^{-q_{D+1} s} }{{\sigma}_{D+1}},\:\:\:\mbox{with}\:\:\: {\sigma}_{D+1} :=e^{-q_{D+1} (t-s_D+D\Delta)}\]
 (and $0$ elsewhere).
Observe that all paths sampled under ${\mathbb Q}$ are necessarily in ${\mathscr Z}(\Delta)$.
The likelihood ratio of such a path reads
\[L=\pi_1 \prod_{i=1}^D\left(\frac{q_{i,i+1}}{q_i}\right)\cdot \left(\prod_{i=1}^{D+1}\sigma_i\right).\]
Performing ${M}_2$ runs, we have thus constructed the estimator, with $L_{{\ell}}$ the likelihood ratio corresponding with the ${{\ell}}$-th sample,
\[\frac{1}{{M}_2}\sum_{{{\ell}}=1}^{{M}_2}  {\mathscr P} (Na, N\phi_t(J_{{\ell}}))L_{{\ell}}.\]
\end{itemize}
As an alternative, one could use the following estimator (in self-evident notation), based on ${M}$ runs under the original and alternative measure:
\[\frac{1}{ {M}}\left(\sum_{{{\ell}}=1}^{ {M}}  {\mathscr P} \left(Na, N\phi_t(J^{({\mathbb P})}_{{\ell}})\right) \,1\{J^{({\mathbb P})}_{{\ell}}\not\in {\mathscr Z}(\Delta)\}+
 {\mathscr P} \left(Na, N\phi_t(J^{({\mathbb Q})}_{{\ell}})\right)L_{{\ell}}
\right).\]

\begin{figure}[t]
\begin{center}
\begin{subfigure}{0.65\linewidth}
\includegraphics[width=\linewidth, height=0.67\linewidth]{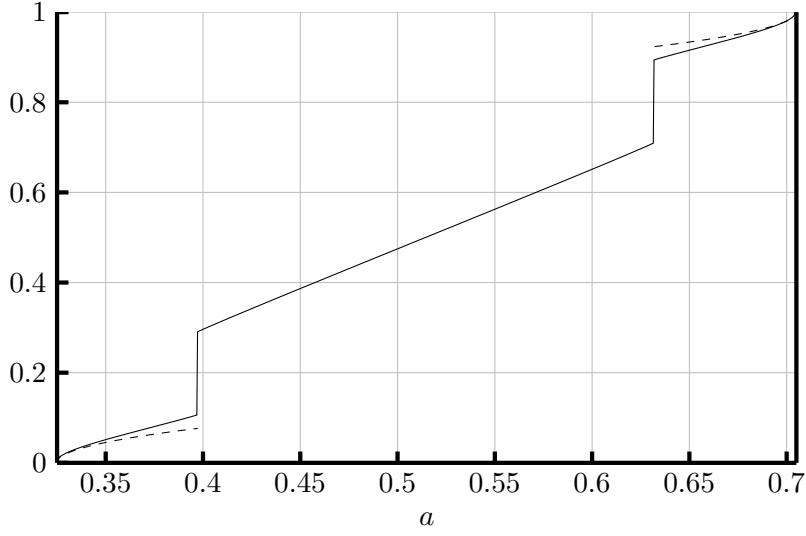}
\end{subfigure}
\caption{\label{figqq}The distribution function ${\mathbb P}(\psi_1(J)\le a)$ for $a\in[a^\amii_1,a^\apii_1]$, dashed the curves $\bar\kappa_1\sqrt{a- a^\amii_1}$ and $1-\bar\kappa_1 \sqrt{a^\apii_1-a}$.}\end{center}
\end{figure}

\begin{figure}[t]
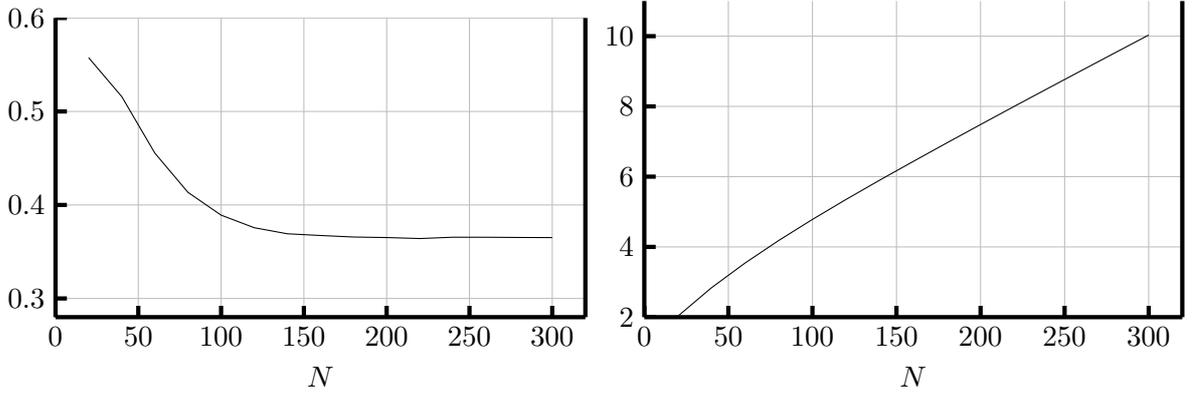

\begin{center}
\begin{subfigure}{0.47\linewidth}
\includegraphics[width=\linewidth, height=0.68\linewidth]{B_fig.tikz}
\end{subfigure}
\begin{subfigure}{0.47\linewidth}
\includegraphics[width=\linewidth, height=0.68\linewidth]{C_fig.tikz}
\end{subfigure}
\caption{\label{figq}  Left panel: $N e^{NI}\,p_1^{(N)}(1) $ for $N\in\{20,40,\ldots,300\}$; right panel: $-\log_{10} p_1^{(N)}(1) $ for $N\in\{20,40,\ldots,300\}$.}\end{center}
\end{figure}

\begin{exmp} \label{EX3}
Following up on Example \ref{EX1}, we consider Model {\sc ii} with $d=2$ and the following choice of the parameters: $\lambda_1=\mu_1=1$, $\lambda_2=2$, $\mu_2=5$, $q_1=q_2=1$, and $t=1$. As it turns out, $s_1=1-\log{\sqrt[4]{2}}$, and
\[a^\apii_1= \int_0^{1-\log {\sqrt[4]{2}}} \lambda_1 e^{-\mu_1 (1-r)}{\rm d}r+\int_{1-\log {\sqrt[4]{2}}}^1  \lambda_2 e^{-\mu_2 (1-r)}{\rm d}r=
\frac{1}{\sqrt[4]{2}}-\frac{1}{e}+\frac{2}{5}\left(1-\left(\frac{1}{\sqrt[4]{2}}\right)^5\right)
,\]
which equals $0.704838$. We focus on the probability $p_1^{(N)}(a)$ that $M^{(N)}(t)$ exceeds $Na$, with $a=1> a^\apii_1$. 
Likewise,
\[a^\amii_1= \int_0^{1-\log {\sqrt[4]{2}}} \lambda_2 e^{-\mu_2 (1-r)}{\rm d}r+\int_{1-\log {\sqrt[4]{2}}}^1  \lambda_1 e^{-\mu_1 (1-r)}{\rm d}r=
\frac{2}{5}\left(\left(\frac{1}{\sqrt[4]{2}}\right)^5-e^{-5}\right)+1 - \frac{1}{\sqrt[4]{2}}
,\]
which equals $0.324588$. Fig.\ \ref{figqq} presents the distribution function of $\psi_1(J)$. Observe that there are atoms
of size $\pi_1 e^{-q_1 t} = (2e)^{-1}\approx 0.183940$ at $1-e^{-1}\approx 0.632120$, and of size $\pi_2 e^{-q_2 t} = (2e)^{-1}\approx 0.183940$ at $\frac{2}{5}(1-e^{-5})\approx 0.397305$; these atoms correspond to the scenarios that the process starts in state 1 (state 2, respectively) and does not leave that state before $t=1$. It is also seen that the shape of ${\mathbb P}(\psi_1(J)\le a^\amii_1+\delta)$ as well as ${\mathbb P}( \psi_1(J)\ge  a^\apii_1-\delta)$ for $\delta$ small is roughly proportional to $\sqrt{\delta}$, in line with results derived earlier in this paper.

By virtue of Thm.\  \ref{st} we know that $N e^{NI}\,p_1^{(N)}(1)$ 
should converge to a constant as $N\to\infty$,
with the decay rate $I$ equal to
 \[I(1\,|\,a^\apii_1)=
-\log a^\apii_1+ a^\apii_1 -1\approx 0.0546252;\]
this convergence is confirmed by the left panel of Fig.\ \ref{figq}. The right panel of Fig.\ \ref{figq} shows the (approximately) exponential decay of $p_1^{(N)}(1) $ (as a function of $N$).
\end{exmp}

\begin{figure}[t]
\begin{center}
\begin{subfigure}{0.6\linewidth}
\includegraphics[width=\linewidth, height=0.56\linewidth]{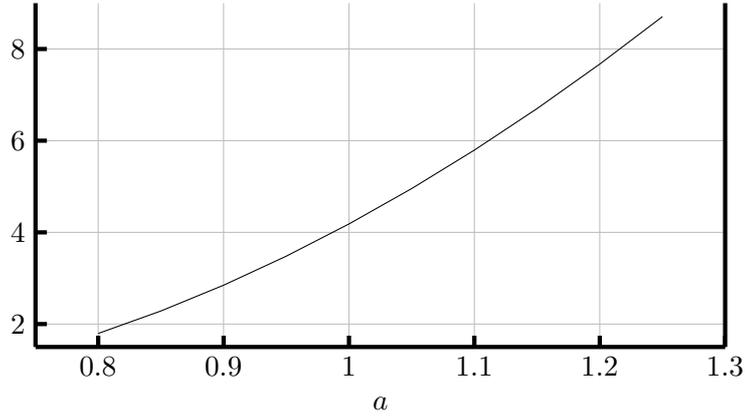}
\end{subfigure}
\caption{\label{figqqq}$-\log_{10} p_1^{(80)}(a)$ for $a\in[0.8,1.25]$.}\end{center}
\end{figure}

\begin{exmp}
In this example we take the same parameters as in Example \ref{EX3}, but fix $N=80$. Our objective is to find, for a given value of $\varepsilon$, the value of $a$ such that $p_1^{(80)}(a)<\varepsilon.$ Then $Na$ could be used as a (somewhat rough) approximation of the number of servers needed in the corresponding finite-server system so as to keep the blocking probability below $\varepsilon$. From Fig.\ \ref{figqqq} we see that e.g.\ for $\varepsilon =10^{-3}$ we need $80\cdot 0.92\approx 74$ servers, and for $\varepsilon =10^{-4}$ we need $80\cdot 0.98\approx 78$ servers.
\end{exmp}

\section{Discussion and concluding remarks}\label{CONCL}
In this paper we have identified the exact asymptotics of the tail distribution of the number of jobs $M^{(N)}(t)$ present in a Markov-modulated infinite-server queue at some time $t>0$; this finding extends earlier obtained logarithmic asymptotics \cite{BKMT,BM}. In the asymptotic regime that we consider, in which  the arrival rates are inflated by a factor $N$, the exact asymptotics are the product of a polynomial function (in $N$) and an exponential function (in $N$). The degree of the polynomial function depends on the number of jumps  the background process makes so as to maximize the (random) Poisson parameter that describes the distribution of $M^{(N)}(t).$

In our paper we have concentrated on the exact asymptotics for the model in which the transition rate matrix $Q$ of the background process is not scaled. A topic for future research could relate to identifying such asymptotics for the setting in which $Q$ is scaled by a factor $N^\alpha.$ For $\alpha=1$ logarithmic asymptotics have been obtained in \cite{DET}, where related results in a more general diffusion setting were derived in \cite{HUANG} building on the framework developed in \cite{LIPTSER}, but these do not seem to lend themselves to a straightforward extension to exact asymptotics. For $\alpha>1$ the system essentially behaves as an ordinary (non-modulated, that is) M/M/$\infty$ queue, and  it is therefore conceivable that its exact asymptotics coincide with those of that M/M/$\infty$ queue.

{\small \section*{Acknowledgment}

The authors would like to thank Mark Peletier and Sorin Pop (both Eindhoven University of Technology) and Marijn Jansen (University of Ghent and University of Amsterdam) for valuable remarks.}

\end{document}